\documentclass[a4paper, 12pt]{article}
\usepackage{enumerate, theorem}
\usepackage{amsmath, amsfonts, amssymb}
\usepackage[height=22.5cm, width=15cm]{geometry}
\usepackage[all]{xy}
\usepackage{mathrsfs, yfonts}

\DeclareMathOperator{\id}{id}

\DeclareMathOperator{\C}{\mathbb{C}}

\newcommand{\A}{\tilde{\mathcal{A}}}
\newcommand{\Ab}{\mathcal{A}}

\newcommand{\parag}[1]{\paragraph{\sc{#1.}}}

\newtheorem{thm}{Theorem}[subsection]
\newtheorem{defn}[thm]{Definition}
\newtheorem{cor}[thm]{Corollary}
\newtheorem{prop}[thm]{Proposition}
\newtheorem{lemma}[thm]{Lemma}

\begin{document}
\title{Algebraic differential equations associated to some polynomials (seconde version)}

\author{Daniel Barlet\footnote{Barlet Daniel, Institut Elie Cartan UMR 7502  \newline
Universit\'e de Lorraine, CNRS, INRIA  et  Institut Universitaire de France, \newline
BP 239 - F - 54506 Vandoeuvre-l\`es-Nancy Cedex.France. \newline
e-mail : daniel.barlet@univ-lorraine.fr}.}

\date{24/02/14.}

\maketitle

\section*{Abstract} 
We compute the Gauss-Manin differential equation for some periods of a polynomial in \ $\C[x_{0},\dots, x_{n}]$ \ with \ $(n+2)$ \ monomials which is not quasi-homogeneous. We show that the corresponding differential equations are defined over \ $\mathbb{Q}[\lambda]$ \ where \ $\lambda \not= 0$ \ is the ``natural'' parameter of such a polynomial. We give also two general  factorization theorems in the algebra \ $\Ab : = \C< z, (\frac{\partial}{\partial z})^{-1}>$ \ for such a differential equation.

\parag{AMS Classification} 32 S 25, 32 S 40.

\tableofcontents

\newpage

\section{Introduction and  statement of the  theorem.}

\subsection{Introduction.}
Since the works of J. Milnor [Mi.68] and E. Brieskorn [Br.70] \ it is a classical problem to compute the Gauss-Manin connection of a holomorphic function. The first examples are the Brieskorn-Pham singularities and then the quasi-homogeneous isolated singularities. We shall give here a simple calculation which shows that in the next case, that is to say  for a non quasi-homogeneous polynomial in \ $n+1$ \ variables with \ $n+2$ \ monomials, it is possible to obtain for any given integral of the type \ $\int_{\gamma_{s}} \ \mu.dx\big/df $, where \ $\mu$ \ is a monomial and \ $\gamma_{s}$ \ a horizontal family of \ $n-$cycles in the fibers of \ $f$, an explicit differential (global) equation. This differential equation is defined over \ $\mathbb{Q}[\lambda]$ \ where \ $\lambda \not= 0$ \ is the parameter of such a polynomial (defined modulo a diagonal change of variables) and has a very special form described in our theorem below. This motivate the study of the kind of differential equation which appears in these examples.  We prove two interesting factorization theorems which apply to such an equation in section 3.

\subsection{Statement of the theorem.}

We consider a polynomial \ $f \in \C[x_{0}, \dots, x_{n}] $ \ which is the sum of \ $n+2$ \ monomials \ $ f = \sum_{j=1}^{n+2} \ m_{j}$ \  where \ $m_{j} : = \sigma_{j}.x^{\alpha_{j}}$, with \ $\sigma_{j}\in \C^{*}$ \ and \ $\alpha_{j}\in \mathbb{N}^{n+1}$. Define the matrix with \ $(n+1)$ \ lines and \ $(n+2)$ \ columns \ $M = (\alpha_{i,j})$ \ and let  \ $\tilde{M}$ \ be the square matrix obtained from \ $M$ \ by adding a first line equal to \ $(1, \dots, 1)$. We shall assume the following condition:
\begin{enumerate}[(C)]
\item The rank of \ $\tilde{M}$ \ is \ $n+2$.
\end{enumerate}

Remark that  this condition is equivalent to the fact  that such a \ $f$ \ is not quasi-homogeneous.

\parag{Example} Let \ $p_{0}, \dots, p_{n}$ \ be integers bigger or equal to  \ $2$ \ and \ $\beta \in (\mathbb{N}^{*})^{n+1}$ \ such that  
$$ \sum_{0}^{n} \frac{\beta_{i}}{p_{i} }\ \not=  1, $$
then
$$ f(x_{0}, \dots, x_{n}) : = \sum_{0}^{n} \ \sigma_{i}.x_{i}^{p_{i}} + \rho.x^{\beta} $$
satisfies the condition \ $(C)$, assuming that \ $\sigma_{0}\dots\sigma_{n}.\rho \not= 0$.$\hfill \square$\\

Consider now  \ $\sigma_{1}, \dots, \sigma_{n+2} \in \C^{*}$ \ and define as before
$$ f(x) : = \sum_{j=1}^{n+2} \sigma_{j}.x^{\alpha_{j}} .$$
A diagonal linear change of variables allows to reduce the study to the case where 
\begin{equation*}
 f(x) = \sum_{j=1}^{n+1} \ x^{\alpha_{j}} \ + \lambda.x^{\alpha_{n+2}} \tag{@}
 \end{equation*}
for some \ $\lambda \in \C^{*}$. So in what follows, we shall assume that \ $m_{j} = x^{\alpha_{j}}$ \ for \ $j \in [1,n+1]$ \ and \ $m_{n+2} = \lambda.x^{n+2}$ \ where \ $\lambda \in \C^{*}$ \ is a  fixed parameter.\\

In our setting it will be convenient to construct differential equations using, instead of the Weyl algebra \ $\C<s, \frac{\partial}{\partial s}>$,  the algebra \ $\Ab : = \C<a,b> \simeq \C<s, \frac{\partial}{\partial s}^{-1}>$ \ defined by the commutation relation \ $a.b - b.a = b^{2}$ \ which is the translation of the Leibnitz rule, as we shall make \ $a$ \ acts on holomorphic functions as multiplication by \ $s$ \ and \ $b$ \ as \ $ (\frac{\partial}{\partial s})^{-1} : = \int_{0}^{s}$. The traduction of our result in usual terms will be immediate (see the remark following the theorem below).\\

For such a polynomial \ $f$ \ we define the unitary sub-algebra \ $A_{f} \subset \C[x_{0}, \dots, x_{n}]$ \ generated by the monomials of \ $f$. It is the image of \ $R : = \C[m_{1}, \dots, m_{n+2}]$ \ by the algebra morphism defined by \ $m_{j} \mapsto x^{\alpha_{j}}, j \in [1,n+1]$ \ and \ $m_{n+2} \mapsto \lambda.x^{\alpha_{n+2}}$, where \ $\lambda\not= 0$ \ is a fixed parameter.
We define also
$$ E_{f} : = \Omega^{n+1}\big/df\wedge d\Omega^{n-1} $$
where \ $\Omega^{k} : =  \C[x_{0}, \dots, x_{n}]\otimes_{\C} \Lambda^{k}((\C^{n+1})^{*})$ \ where \ $(\C^{n+1})^{*} : = \oplus_{i=0}^{n} \ \C.dx_{i}$. \\
Then \ $E_{f}$ \  is a left  \ $\Ab-$module where \ $a$ \ acts by multiplication by \ $f$ \ and where \ $b$ \ acts as \   $df\wedge d^{-1}$.We shall fix in the sequel a monomial \ $\mu : = x^{\beta}$ \ in \ $ \C[x_{0}, \dots, x_{n}]$ \ and  we shall denote \ $\varphi_{\mu} : R \to  E_{f}$ \ the map given by \ $ m_{j} \mapsto [m_{j}.\mu.dx]$.

\begin{thm}\label{base}
In the situation above there exist  positive integers \ $d$ \ and \ $h$ \ such that for any choice of the monomial \ $\mu \in \C[x_{0}, \dots, x_{n}] $ \ we may construct in a natural way (explained below) homogeneous elements in \ $(a,b)$ \ of degree \ $d$ \ and \ $d+h$ \ in the algebra \ $\Ab$, noted \ $P_{d}(\mu)$ \ and \ $P_{d+h}(\mu)$,  which are monic in \ $a$,  and such that the element
$$ P(\mu) : = P_{d+h}(\mu) - c.\lambda^{r}.P_{d}(\mu) \in \Ab $$
with some  \ $r \in \mathbb{Z}^{*}$ \ and \ $c \in \mathbb{Q}^{*}$, has the following properties :
\begin{enumerate}
\item Assume that on a connected open set \ $S$ \ in \ $\C^{*}$ \ with the origin on its boundary,  we have an open set \ $U \subset f^{-1}(S)$ \ such that \ $f : U \to S$ \ is a \ $\mathscr{C}^{\infty}-$fibration with fiber \ $F$;  then for any  compact \ $n-$cycle \ $\gamma \in H_{n}(F, \C)$, the holomorphic function on \ $S$ \ defined by (see [M.74])
$$ F_{\mu,\gamma}(s) : = \int_{\gamma_{s}} \mu.dx\big/df  $$
where \ $\gamma_{s}$ \ is the horizontal family of cycles in the fibers of \ $f$ \  associated to \ $\gamma$, satisfies the global algebraic differential equation 
$$ P(\mu)[F_{\mu,\gamma}] = 0 .$$
\item Each of the elements \ $P_{d}$ \ and \ $P_{d+h}$ \ are products in the algebra \ $\Ab$ \ of elements of the form \ $a - r.b$ \ where \ $r$ \ is a rational number.
\item The elements \ $P_{d}$ \ and \ $P_{d+h}$ \ are computable by an easy linear algebra algorithm.
\end{enumerate}
\end{thm}

\parag{Remark} To reach the algebraic differential equation in the usual way just look at \ $b^{-d-h}P(\mu) $ \ which is the sum of a monic polynomial of degree \ $d+h$ \  (with rational roots) in \ $s.\frac{\partial}{\partial s}$, (see the property 2.)  minus \ $c.\lambda^{r}.(\frac{\partial}{\partial s})^{h}$ \ composed with  a monic polynomial of degree \ $d$ \  in \ $s.\frac{\partial}{\partial s}$ \ (also with rational roots). The order of this differential equation is \ $d+h$ \ and  the coefficient of \ $(\frac{\partial}{\partial s})^{d+h}$ \ is equal to \ $s^{d+h} - c.\lambda^{r}.s^{d}$, showing that the singular points of this differential equation  outside the origine  are  solutions of \ $s^{h} - c.\lambda^{r} = 0$. $\hfill \square$\\

%We shall see in the proof that the construction of \ $P_{d}(\mu)$ \ and \ $P_{d+h}(\mu)$ \ is given by a simple algorithm of linear algebra.

\parag{Example} Assume that \ $p \in f^{-1}(0)$ \ is a singular point of  a polynomial \ $f$ \ satisfying the condition \ $(C)$. Then choose for \ $U$ \  a Milnor ball for \ $f$ \ at \ $p$ \ and for \ $S$ \ the corresponding punctured disc with center \ $0$. Then they fulfill the hypothesis in the theorem, and so we obtain, for each given monomial \ $\mu$, an algebraic global differential equation satisfied by the function \ $F_{\mu,\gamma}$ \ for any choice of \ $\gamma$ \ in the \ $n-$th homology group of the Milnor fiber of \ $f$ \ at the point \ $p$.$\hfill \square$\\

\section{Proof of the theorem.}

\subsection{Some \ $\Ab-$modules.}

We consider a polynomial \ $f \in \C[x_{0},\dots,x_{n}]$ \ with \ $n+2$ \ monomials \ $m_{1}, \dots, m_{n+2}$.  So \ $ f = \sum_{j=1}^{n+2} \ m_{j}$. We write \ $m_{j} : = \sigma_{j}.x^{\alpha_{j}}$ \ with \ $\sigma_{j} \in \C^{*}$, and \  $\alpha_{j}\in \mathbb{N}^{n+1}$. Using a linear diagonal change of variables and reordering the variables if necessary we shall assume the following conditions :
\begin{enumerate}[i)]
\item We have \ $\sigma_{j} = 1$ \ for all \ $j \in [1,n+1]$ \ and \ $\sigma_{n+2} = \lambda \in \C^{*}$.
\item We assume that \ $\alpha_{1}, \dots, \alpha_{n+1}$ \ are linearly independent in \ $\mathbb{Q}^{n+1}$.
\end{enumerate}
Then we shall write
\begin{equation*}
\alpha_{n+2} = \sum_{j=1}^{n+1} \ \rho_{j}.\alpha_{j} \quad {\rm with} \quad \rho_{j} \in \mathbb{Q}. \tag{@}
\end{equation*}
We define
\begin{align*}
& H : = \{j \in [1,n+1] \ / \ \rho_{j} = 0 \} ;\\
& J_{+} : = \{j \in [1,n+1] \ / \ \rho_{j} > 0 \} ;\\
& J_{-} : = \{j \in [1,n+1] \ / \ \rho_{j} < 0 \}.
\end{align*}
Let \ $\vert r\vert$ be the smallest positive integer such that \ $\vert r\vert.\rho_{j} : = p_{j}$ \ is an integer for each \ $j \in [1,n+1]$. Write now the relation above as
$$ \vert r\vert.\alpha_{n+2} + \sum_{j \in J_{-}} (-p_{j}).\alpha_{j} = \sum_{j \in J_{+}} p_{j}.\alpha_{j}.$$
Now define \ $ d+h$ \ and \ $d$ \ as respectively the supremum and infimum of the two positive integers  \ $\vert r\vert + \sum_{j \in J_{-}} \ (-p_{j})$ \ and \ $  \sum_{j \in J_{+}}\  p_{j} $. Here \ $d$ \ and \ $h$ \ are positive because  the positivity of \ $d$ \ is consequence of the fact that \ $\vert r\vert \geq 1$ \ and that at least one\ $p_{j}$ \ is positive.\\
The positivity of \ $h$ \ is consequence of the fact that the equality of these two integers would  imply that the first line in \ $\tilde{M}$ \ satisfies the same linear relation \ $(@)$ \ than all the other lines in \ $\tilde{M}$, contradicting our hypothesis \ $(C)$.\\
Now the relation above gives the relation between the monomials \ $(m_{j})_{j\in [1,n+2]}$ :
$$ m_{n+2}^{\vert r\vert}.\prod_{j \in J_{-}} m_{j}^{-p_{j}} = \lambda^{\vert r\vert}.\prod_{j \in J_{+}} m_{j}^{p_{j}} $$
and we shall write it 
 \begin{equation*}
  m^{\Delta} = \lambda^{r}.m^{\delta} \tag{@@}
  \end{equation*}
 where \ $\Delta$ \ and \ $\delta$ \ are in \ $\mathbb{N}^{n+2}$ \ of respective weight \ $d+h$ \ and \ $d$, and with zero component for each \ $h \in H$. Note that we have \ $r = \pm \vert r\vert $ \ and so that\ $r $ \ is in \ $\mathbb{Z}^{*}$.\\
  
We shall also use the following observation later on :
\parag{Observation} For \ $j \in [1,n+1]$ \ the \ $j-$th element of the first column of the matrix \ $\tilde{M}^{-1}$ \ is zero if and only if \ $j$ \ is in \ $H$.\\

Consider the  free polynomial algebra \ $R : = \oplus_{q=0}^{\infty}  \ R_{q} = \C[m_{1}, \dots, m_{n+2}]$ \ graded by the total degree in the variables \ $m_{1}, \dots, m_{n+2}$ \ and consider  the algebra morphism  \ $\varphi : R \to \C[x_{0},\dots,x_{n}]$ \ sending the variable \ $m_{j}$ \ in \ $R$ \ to the monomial \ $x^{\alpha_{j}}$ \ in \ $\C[x_{0},\dots,x_{n}]$ \ for \ $j \in [1,n+1]$ \ and \ $m_{n+2}$ \ to \ $\lambda.x^{\alpha_{n+2}}$ \ where \ $\lambda \in \C^{*}$ \ is a fixed parameter.\\
We shall need an extra variable \ $\mu$ \ and we associated to it a multi-index \ $\beta \in \mathbb{N}^{n+1}$. Then we consider the free, rank \ $1$ \ $R-$module \ $R.\mu$ \ and the \ $R-$linear map \ $\varphi_{\mu} : R.\mu \to \C[x_{0},\dots,x_{n}]$ \ by sending \ $\mu$ \ to \ $x^{\beta}$.\\
We construct now on the two vector spaces  \ $R.\mu$ \ and \ $ \C[x_{0},\dots,x_{n}] $ \  the  \ $\C-$linear operators\ $b_{i}, i \in [0,n]$ \ and \ $a$ \ as follows :
\begin{enumerate}[i)]
\item For \ $\gamma \in \mathbb{N}^{n+1}$ \ and \ $i \in [0,n]$ \ let \ $\Gamma_{i}(\gamma,\mu) : = \beta_{i} +  \sum_{j=1}^{n+k} \alpha_{i,j}.\gamma_{j} $ \ and define  in the vector space  \ $R.\mu$
$$(\Gamma_{i}(\gamma,\mu)+1).b_{i}(m^{\gamma}.\mu) = \sum_{j=1}^{n+2} \alpha_{i,j}.m_{j}.m^{\gamma}.\mu \quad\quad  {\rm and} \quad \quad a(m^{\gamma}.\mu) =  \sum_{j=1}^{n+2} m_{j}.m^{\gamma}.\mu$$
\item For \ $\zeta \in \mathbb{N}^{n+1}$ \ and \ $i \in [0,n]$ \  define on the vector space \ $\C[x_{0},\dots,x_{n}]$ :
$$ (\zeta_{i} +1).b_{i}(x^{\zeta}) = \sum_{j=1}^{n+1} \alpha_{i,j}.x^{\alpha_{j}+ \zeta} \  + \  \lambda.x^{\alpha_{n+2}+\zeta} = x^{\zeta}.(x_{i}.\frac{\partial f}{\partial x_{i}}) \quad \quad {\rm and} \quad \quad a(x^{\zeta}) : = f.x^{\zeta}$$
\end{enumerate}
Note that  the operator \ $b_{i}$ \ on \ $\C[x_{0},\dots,x_{n}]$  \ may also be written as :
$$ P \mapsto \frac{\partial f}{\partial x_{i}}.\int_{0}^{x_{i}} P.dt = \left(x_{i}.\frac{\partial f}{\partial x_{i}}\right).\frac{1}{x_{i}}.\int_{0}^{x_{i}} P.dt  = \left(\sum_{j=1}^{n+2} \alpha_{i,j}.m_{j}\right).\big[\frac{1}{x_{i}}.\int_{0}^{x_{i}} P.dt \big].$$
Remark that \ $\C[x_{0},\dots,x_{n}].x^{\beta} \subset  \C[x_{0},\dots,x_{n}] $ \  is stable by the operators \ $b_{i}$ \ and \ $a$ \ for any choice of \ $\beta$.\\
The commutation relations \ $a.b_{i} - b_{i}.a = b_{i}^{2}$ \ are easy to verify on \ $R.\mu$ \ and \ $\C[x_{0},\dots,x_{n}]$. It is also easy to see that \ $\varphi_{\mu}$ \    commutes with these operators.\\
We shall denote \ $V_{p}(\mu)$ \ the vector subspace of \ $R_{p}.\mu$ \ defined inductively as follows:
\begin{enumerate}
\item \ $V_{0}(\mu) = \{0\}$.
\item \ $V_{p+1}(\mu) : = \sum_{i\not= 0} \ (b_{i} - b_{0}).R_{p}.\mu +  b_{0}.V_{p}(\mu)$.
\end{enumerate}
Remark that, by definition, we have  $$b_{i}.V_{p}(\mu) \subset  (b_{i} - b_{0}).V_{p}(\mu) + b_{0}.V_{p }(\mu) \subset V_{p+1}(\mu)$$
 for all \ $i \in [0,n]$ \ and for all \ $p \geq 0$. We have also \ $a.V_{p}(\mu) \subset V_{p+1}(\mu)$ \ because of the identity 
$$ a.(b_{i}-b_{0}) = (b_{i}- b_{0}).a + (b_{i}- b_{0}).b_{i} + b_{0}.(b_{i}- b_{0}) .$$
Then we define the graded vector space \ $S(\mu) : = \oplus_{p \geq 0} \ S_{p}(\mu)$ \ where we put 
  $$S_{p}(\mu) : = R_{p}.\mu\big/V_{p}(\mu).$$
  The quotient map \ $\pi_{\mu} : R.\mu \to S(\mu)$ \ is graded.\\
  
  Note that the algebra \ $\Ab$ \  has a natural grading
$$ \Ab = \oplus_{q=0}^{\infty} \Ab_{q} $$
where \ $\Ab_{q}$ is the vector subspace of homogeneous elements in \ $a,b$ \ of degree \ $q$.

\parag{Remark} An element \ $P \in \Ab_{q}$ \ which is monic in \ $a$ \ can be written
 $$P = (a - r_{1}.b)\dots (a - r_{q}.b)$$
 for some suitable choice of complex numbers \ $r_{1}, \dots, r_{q}$ ; see [B.09] \ or write \ $b^{-q}.P$ \ as a polynomial in \ $b^{-1}.a$.\hfill $\square$\\

On \ $S(\mu)$ \ we have a structure of left  graded \ $\Ab-$module deduced from the action of the operators \ $a$ \ and \ $b : = b_{0} = \dots = b_{n}$ \  which are equal on \ $S(\mu)$. Note that if we consider \ $R.\mu$ \ as a left \ $\Ab-$module with the actions of \ $a$ \ and \ $b_{0}$, the quotient map \ $\pi_{\mu} : R.\mu \to S(\mu)$ \ becomes a surjective map of  graded left \ $\Ab-$modules.

\begin{prop}\label{new}
Fix \ $\beta \in \mathbb{N}^{n+1}$ \ and  \ $q \in \mathbb{N}$. For each \ $\gamma \in \mathbb{N}^{n+2}$ \ such that  \ $\vert \gamma \vert = q$, there exists an element \ $\chi_{q}(\mu)(m^{\gamma}) \in \Ab_{q}$, which is a product of homogeneous elements in (a,b) of degree \ $1$ \ of the form \ $\eta.a + \theta.b$ \ with \ $\eta $ \ and \ $\theta$ \ in  \ $ \mathbb{Q}$, such that we have in \ $S(\mu)$ \ the equality :
$$ \pi_{\mu}(m^{\gamma}.\mu) = \chi_{q}(\mu)(m^{\gamma}).\mu .$$
Moreover, for  every\ $\gamma$ \ such that \ $\gamma_{h} = 0$ \ for each \ $h \in H $, the corresponding element \ $\chi_{q}(\mu)(m^{\gamma}) \in \Ab_{q}$ \ is monic in \ $a$ \ up to a non zero rational number.
\end{prop}

\parag{Proof} Consider the following equalities in \ $S(\mu)$ :

\begin{align*}
& a[\pi_{\mu}(m^{\gamma}.\mu)] = f.m^{\gamma}.\mu.dx = \sum_{j=1}^{n+2} \ \pi_{\mu}(m_{j}.m^{\gamma}.\mu) \\
& (\Gamma_{i}(\gamma,\mu)+1).b[\pi_{\mu}(m^{\gamma}.\mu)] = \sum_{j=1}^{n+2} \ \alpha_{i,j}.\pi_{\mu}(m_{j}.m^{\gamma}.\mu) \quad \forall i \in [0,n]. \tag{@@@}
\end{align*}
as a linear system with matrix \ $\tilde{M}$ \ and unknown the  column vector 
 $$\big( \pi_{\mu}(m_{j}.m^{\gamma}.\mu) \big), j \in [1,n+2]$$ 
 in \ $S(\mu)^{n+2}$. The invertibility of \ $\tilde{M}$ \ (on the field \ $\mathbb{Q}$) gives that for each \ $j \in [1,n+2]$ \ we may write \ $\pi_{\mu}(m_{j}.m^{\gamma}.\mu) = (\eta_{j}.a + \theta_{j}.b).\pi_{\mu}(m^{\gamma}.\mu)$ \ for some rational numbers \ $\eta_{j}$ \ and \ $\theta_{j}$. Assuming that, for each \ $\gamma$ \ with \ $\vert \gamma\vert = q$ \ the element \ $\chi_{q}(\mu)(m^{\gamma})$ \ is a product of \ $q$ \ homogeneous elements in (a,b) of degree \ $1$ \ of the form \ $\eta.a + \theta.b$ \ with \ $\eta$ \ and \ $\theta$ \ in \ $ \mathbb{Q}$, we conclude  that for each \ $\tilde{\gamma}$ \ with \ $\vert\tilde{\gamma}\vert = q+1$ \ the element  \ $\chi_{q+1}(\mu)(m^{\tilde{\gamma}})$ \ also satisfies this condition.\\
Moreover, as the first column of \ $\tilde{M}^{-1}$ \ has a non zero \ $j-$th component for \ $j \not\in H$ \ (see the observation above),  we have \ $\eta_{j} \not= 0$ \ for each \ $j \not\in H$. So, assuming inductively that, for any \ $\gamma$ \ such that \ $\vert \gamma \vert = q$, with zero \ $j-$th component for \ $j \in H$, we have already constructed  a suitable \ $ \chi_{q}(\mu)(m^{\gamma}) \in \Ab_{q}$ \  this implies that for any \ $\tilde{\gamma} \in \mathbb{N}^{n+2}$ \ such that \ $\vert \tilde{\gamma} \vert = q+1$ \ and with no component in \ $H$, we may write \ $m^{\tilde{\gamma}} = m_{j}.m^{\gamma}$ \ for some \ $ j \not\in  H$ \ and some \ $\gamma \in \mathbb{N}^{n+2}$ \ such that \ $\vert \gamma \vert = q$ \ has zero components in \ $H$, and define \ $\chi_{q+1}(\mu)(m^{\tilde{\gamma}}) : = (\eta_{j}.a + \theta_{j}.b).\chi_{q}(\mu)(m^{\gamma})$ \ where \ $\eta_{j}$ \ is in \ $\mathbb{Q}^{*}$. So we conclude that for such a \ $\gamma$ we may choose \ $\chi_{q+1}(\mu)(m^{\tilde{\gamma}}) $ \ as a monic element in \ $a$ \ belonging to \ $\Ab_{q+1}$, up to a non zero rational number.$\hfill \blacksquare$\\

Now using only the actions of \ $a$ \ and \ $b_{0}$ \ gives a structure of left \ $\Ab-$module on  \ $R.\mu$ \ and a graded \ $\Ab-$linear map \ $.\mu : \Ab \to R.\mu$ \ sending \ $1$ \ to \ $\mu$. Composed with the quotient map \ $\pi_{\mu}$ \ we obtain a graded \ $\Ab-$linear map \ $\pi_{0} : \Ab \to S_{\mu}$. The following corollary is an obvious consequence of the proposition above.

\begin{cor}\label{surj.}
The graded\ $\Ab-$linear map \ $\pi_{0} : \Ab \to S_{\mu}$ \ is surjective.
\end{cor}

\subsection{End of the proof of the theorem \ref{base}}

The following lemma is now easy.

\begin{lemma}\label{ quot.}
The map \ $\varphi_{\mu}$ \ induces a commutative diagram \\
\begin{equation*} 
\xymatrix{\Ab \ar[rd]_{\pi_{0}} \ar[r]^{.\mu} & R.\mu \ar[d]_{\pi} \ar[r>]^{\varphi_{\mu}} &\C[x_{0},\dots,x_{n}]  \ar[d]_{\pi'}\\
\quad &S_{\mu} \ar[r]^{\psi_{\mu}} & E_{f}} \\
\end{equation*}
where the map \ $\psi_{\mu}$ \ is \ $\Ab-$linear. The image of \ $\psi_{\mu}$ \ is equal to   \ $\Ab.[x^{\beta}.dx] \subset E_{f}$.
\end{lemma}

\parag{proof} The \ $\Ab-$linearity is a direct consequence of the compatibility of the map \ $\varphi_{\mu}$ \ with the operators \ $a$ \ and \ $b_{i}, i \in [0,n]$ \ and the fact that these operators induce \ $a$ \ and \  $b$ \ on \ $S(\mu)$ \ and on \ $E_{f}$. The identification  of the image of \ $\psi_{\mu}$ \ is then  a consequence of the corollary of the proposition \ref{new} and of the commutativity of the diagram above. $\hfill \blacksquare$\\

As the relation \ $(@@)$ \ holds in \ $\C[x_{0}, \dots, x_{n}]$, we have also \ $m^{\Delta}.x^{\beta}.dx = \lambda^{r}.m^{\delta}.x^{\beta}.dx$ \ in \ $E_{f}$. As \ $\Delta$ \ and \ $\delta$ \ have no components in \ $H$, there exists elements \ $P_{d+h}$ \ and \ $P_{d}$ \ respectively in \ $\Ab_{d+h}$ \ and \ $\Ab_{d}$ \ which are equal, modulo some invertible rational numbers, to \ $\chi_{d+h}(\mu)[m^{\Delta}]$ \ and \ $\chi_{d}(\mu)[m^{\delta}]$. Then we obtain in \ $E_{f}$ \ the equality
$$ \big(P_{d+h} - c\lambda^{r}.P_{d}\big)[x^{\beta}.dx] = 0 $$
for some non zero rational number \ $c$, thanks to the \ $\Ab-$linearity of the map \ $\psi_{\mu}$.

Now assume the existence of \ $S$ \ and \ $U$ \ as in the property 1. of the statement of the theorem. Note that we may assume \ $S$ \ simply connected.  Then define \ $\mathcal{O}_{b}(S)$ \ as the vector space of   holomorphic functions on \ $S$ \ such that \ $\vert s\vert^{1+\varepsilon}.\vert f(s)\vert \to 0$ \ when \ $s \to 0$, for some \ $\varepsilon > 0$.  We want to show that the \ $\C-$linear map
$$ I :  E_{f} \longrightarrow \mathcal{O}_{b}(S) \otimes_{\C}H_{n}(F,\C)^{*} $$
   given by \ $\mu \mapsto \big(\gamma \mapsto \int_{\gamma_{s}} \ \frac{\mu.dx}{df} \big) $ \ is \ $\Ab-$linear. Here the actions on \ $\mathcal{O}_{b}(S) \otimes_{\C}H_{n}(F,\C)^{*}$ \ are given by the rules \ $a(g\otimes v) : = s.g\otimes v$ \ and \ $b(g\otimes v) : = (\int_{0}^{s} g)\otimes v$. \\
   The commutation with \ $a$ \ of the map \ $I$ \  is clear. \\
    The  commutation of  the map \ $I$ \ with \ $b$ \   is consequence of the derivation formula for the integral of a holomorphic \ $n-$form :
$$ \frac{\partial}{\partial s}\big(\int_{\gamma_{s}} \ \xi \big) = \int_{\gamma_{s}} \frac{d\xi}{df}.$$
Note that for a holomorphic \ $n-$form the function \ $\int_{\gamma_{s}} \ \xi$ \ has only bounded terms in its asymptotic expansion when \ $s \to 0$ \ (see [M.75]), so its derivative is in the space \ $\mathcal{O}_{b}(S)$ \ defined above because such a function admits a ``standard'' convergent asymptotic (multivalued) expansion when \ $s \to 0$. $\hfill \blacksquare$

\subsection{Integral dependance of \ $f$.}

\begin{lemma}
In the situation of the theorem \ref{base} the element \ $f$ \ is integrally dependent on the subring \ $\C[x_{0}.\frac{\partial f}{\partial x_{0}}, \dots, x_{n}.\frac{\partial f}{\partial x_{n}}] \subset \C[x_{0}, \dots, x_{n}]$.
\end{lemma}

\parag{proof} Let \ $F$ \ be the vector with components \\
 $f, x_{0}.\frac{\partial f}{\partial x_{0}}, \dots, x_{n}.\frac{\partial f}{\partial x_{n}}$ \ and \ $W$ \ the vector with components \ $m_{1}, \dots, m_{n+2}$; then we have
$$ F = \tilde{M}.W $$
and this allows to compute each \ $m_{j}$ \ as a linear combination of the components of \ $F$ \ as \ $\tilde{M}$ \ is invertible. Moreover, for \ $j \not\in H$, the coefficient of \ $f$ \ in this linear combination does not vanishes. Replacing in the equation \ $(@@)$ : 
$$ m^{\Delta} = c.\lambda^{r}.m^{\delta}$$
gives a monic degree \ $d+ h = \vert \Delta\vert$ \ polynomial in \ $f$ \ with coefficients in the ring \ $\C[x_{0}.\frac{\partial f}{\partial x_{0}}, \dots, x_{n}.\frac{\partial f}{\partial x_{n}}]$ \ which vanishes identically. $\hfill \blacksquare$\\

For instance, let \ $f : = x^{2} + y^{3} + z^{4} + \lambda.x.y.z $. Then let \ $X : = f - \frac{1}{2}x.f'_{x} -\frac{1}{3}y.f'_{y} - \frac{1}{4}z.f'_{z}$. The computation gives
$$  2^{6}.3^{4}.4^{3}.X^{12} = \lambda^{12}.(X + 13x.f'_{x})^{6}(X + 13y.f'_{y})^{4}(X + 13z.f'_{z})^{3} $$
which gives  an integral dependance  relation of degree \ $13$ \ of \ $f$ \ on the ring \ $\mathbb{Q}[x.f'_{x}, y.f'_{y}, z.f'_{z}]$, for \ $\lambda$ \ fixed in \ $\mathbb{Q}^{*}$.

\parag{Remark} In fact we have shown that,  in general, for \ $\lambda$ \ given in \ $\mathbb{Q}^{*}$,  $f$ \ is integrally dependant on the ring \ $\mathbb{Q}[x_{0}.\frac{\partial f}{\partial x_{0}}, \dots, x_{n}.\frac{\partial f}{\partial x_{n}}]$. $\hfill \square$

\subsection{Examples.}

First let me give an example with \ $3$ \ variables, small (total) degree \ $ = 6$ \  such that \ $d = 61$ \ and \ $h = 15$.\\
Let \ $f = x^{3}.y + y^{4}.z + z^{5}.x + \lambda.x^{2}.y^{2}.z^{2}$. Let \ $\alpha_{j}, j\in [1,4]$, the exponents of the monomials of \ $f$. Note that the singularity of \ $f$ \ at the origin is not isolated.\\
 The equality
$$ 61.\alpha_{4} =  34.\alpha_{1} + 22.\alpha_{2} + 20.\alpha_{3} $$
gives the relation
\begin{equation*}
  (\lambda.x^{2}.y^{2}.z^{2})^{61} = \lambda^{61}.(x^{3}.y )^{34}.(y^{4}.z)^{22}.(z^{5}.x)^{20} . \tag{R}
\end{equation*}
Of course to compute explicitely \ $P : P_{76} - c.\lambda^{-61}.P_{61}$ \ is a little tedious and without real interest. But computing the monomials in the relation \ $(R)$ \ modulo \ $b.\Ab$ \ gives easily the value of the rational number \ $c$. And so  the non zero critical values of \ $f$ \ are solutions of the equation 
 $$s^{15} = c.\lambda^{-61} \quad\quad \rm{ with} \quad\quad c = -\frac{(61)^{61}.(15)^{15}}{(34)^{34}.(22)^{22}.(20)^{20}}.$$ 
This  is not simple to see by a direct computation.\\
  
%Note that if \ $\zeta^{61} = 1$ \ the change of variable 
 %$$x \mapsto \zeta^{16}.x, \qquad y \mapsto \zeta^{13}.y,\qquad  z \mapsto \zeta^{9}.z$$
% sends \ $f_{\lambda}$ \ to \ $f_{\zeta^{15}.\lambda}, \ dx\wedge dy\wedge dz \mapsto \zeta^{38}.dx\wedge dy\wedge dz $ \ and let invariant  the differential equation associated to the class \ $[dx\wedge dy\wedge dz].$\\

We shall give now a family of examples with \ $n = 3$ \ (so \ $4$ \ variables) with some symetry. Let \ $x_{0},x_{1},x_{2},x_{3}$ \ be the coordinates in \ $\C^{4}$, and note \ $\sigma$ \ the circular permutation of the coordinates defined by \ $\sigma(x_{i}) = x_{i+1}, i\in \mathbb{Z}\big/4.\mathbb{Z}$. Now consider a multi-index \ $\alpha \in \mathbb{N}^{4}$ \ with the following conditions :
\begin{enumerate}[i)]
\item \ $\vert \alpha\vert \geq 5$.
\item The matrix \ $M_{0}$ \ given by \ $M_{0}(i,j) : = \alpha_{i+j}$ \ with \ $(i,j) \in  \mathbb{Z}\big/4.\mathbb{Z}$ \ has rank \ $4$\footnote{Note that this is true if and only if \ $\alpha_{0} + \alpha_{2} \not= \alpha_{1} + \alpha_{3}$.}.
\end{enumerate}
Then we shall consider the polynomial
$$ f(x) : = \sum_{j=0}^{3} \ \sigma^{j}(x^{\alpha}) \ + \lambda.x^{\bf{1}} $$
where \ $\lambda$ \ is a non zero complex parameter and where  we use the following conventions : for a monomial \ $\mu$ \ in the variables  \ $x_{0},x_{1},x_{2},x_{3}$ \ we define  \ $\sigma(\mu)$ \ to be the monomial obtained from \ $ \mu$ \ by applying the circular permutation \ $\sigma$ \ to the variables, and  we note \ $\bf{1}$ \ the multi-index \ $(1,1,1,1)$.\\

The relation between the five monomials in \ $f$ \ is simply 
 $$\lambda^{\vert\alpha\vert}.\prod_{j=0}^{3} \sigma^{j}(x^{\alpha}) = (\lambda.x^{\bf{1}})^{\vert\alpha\vert}.$$
 So we will compute the images of the monomials of both sides of this relation by the map \ $\psi_{\bf{1}}$ \ in \ $E_{f}$ \ as \ $P_{d+h}[dx]$ \ and \ $P_{d}[dx]$ \ where \ $P_{m}$ \ are in \ $\Ab_{m}$ \ and monic in \ $a$ \ up to a non zero constant. Here we have \ $d = 4$ \ and \ $h = \vert \alpha\vert - 4$.\\

\parag{Computation of \ $\prod_{j=0}^{3} \sigma^{j}x^{\alpha}$} It will be convenient to use the notation \ $B(x^{\gamma})$ \ for the vector with components \ $(\gamma_{i} +1).b_{i}(x^{\gamma}), i \in [0,3] \simeq  \mathbb{Z}\big/4.\mathbb{Z}$, and \ $\Sigma(x^{\gamma})$ \ for the vector with components \ $\sigma^{j}(x^{\gamma}), j \in [0,3] \simeq  \mathbb{Z}\big/4.\mathbb{Z}$,  when \ $x^{\gamma}$ \ is a monomial. Then we have, for each monomial \ $m \in \C[x_{0}, \dots, x_{n}]$ \ the equation
$$ B(m) = m.\big[M_{0}.\Sigma(x^{\alpha}) + \lambda.x^{\bf{1}}.v_{0} \big]  $$
where \ $v_{0}$ \ is the column  vector \ $v_{0} : = (1, 1, 1, 1)$.
We shall use successively the values 
 $$1, \qquad x^{\alpha},\qquad x^{\alpha}.\sigma(x^{\alpha}), \qquad x^{\alpha}.\sigma(x^{\alpha}).\sigma^{2}(x^{\alpha})$$
  for \ $m$ \ in this formula. But remark that, because we have \ $b_{i}(m) = b(m)$ \ for each \ $i \in [0,3]$, we may write in \ $S(1) = \C[m_{1}, \dots, m_{5}]\big/V(1)$:
\begin{align*}
& B(1) = b(1).v_{0} , \quad B(x^{\alpha})  = b(x^{\alpha}).(v_{0} + M_{0}.v_{1}) \\
&  B(x^{\alpha}.\sigma(x^{\alpha})) = b(x^{\alpha}.\sigma(x^{\alpha})).(v_{0} + M_{0}.v_{2}) +  \lambda.x^{\bf{1}}.x^{\alpha}.\sigma(x^{\alpha}).v_{0}\\ 
&   B(x^{\alpha}.\sigma(x^{\alpha}).\sigma^{2}(x^{\alpha})) = b(x^{\alpha}.\sigma(x^{\alpha}).\sigma^{2}(x^{\alpha})).(v_{0} + M_{0}.v_{3}) + \lambda.x^{\bf{1}}.x^{\alpha}.\sigma(x^{\alpha}).\sigma^{2}(x^{\alpha}).v_{0}  \\
\end{align*}
with the following values for the (column) vectors 
 $$ v_{0} : = \ ^{t}(1, 1, 1, 1), \quad  v_{1} : = \ ^{t}(1, 0, 0, 0),  \quad v_{2}  = \ ^{t}(1, 1, 0, 0)  \quad {\rm and} \quad  v_{3} = \ ^{t}(1, 1, 1, 0) .$$
 Now remark that \ $ M_{0}.v_{0} = \vert\alpha\vert.v_{0}$ \ and recall that \ $M_{0}$ \ is invertible. This allows to deduce the following equalities in \ $E_{f} $:
 \begin{align*}
 & \vert\alpha\vert.\Sigma(x^{\alpha})) = (b(1) -  \lambda.x^{\bf{1}}).v_{0} \\
 & \vert\alpha\vert.x^{\alpha}.\Sigma(x^{\alpha}) = (b(x^{\alpha}) -  \lambda.x^{\bf{1}}.x^{\alpha}).v_{0} + \vert\alpha\vert.b(x^{\alpha}).v_{1} \\
 &  \vert\alpha\vert.x^{\alpha}.\sigma(x^{\alpha}).\Sigma(x^{\alpha}) = \Big(b(x^{\alpha}.\sigma(x^{\alpha})) -  \lambda.x^{\bf{1}}.x^{\alpha}.\sigma(x^{\alpha})\Big).v_{0} + \vert\alpha\vert.b(x^{\alpha}.\sigma(x^{\alpha})).v_{2} \\
 & \vert\alpha\vert.x^{\alpha}.\sigma(x^{\alpha}).\sigma^{2}(x^{\alpha}).\Sigma(x^{\alpha}) = \Big(b(x^{\alpha}.\sigma(x^{\alpha}.\sigma^{2}(x^{\alpha})) -  \lambda.x^{\bf{1}}.x^{\alpha}.\sigma(x^{\alpha}.\sigma^{2}(x^{\alpha})\Big).v_{0}  + \\
 & \qquad\qquad    \qquad\qquad  +  \vert\alpha\vert.b(x^{\alpha}.\sigma(x^{\alpha})).\sigma^{2}(x^{\alpha}).v_{3} 
\end{align*} 
The first equality implies the equality of \ $\sigma^{j}(x^{\alpha})$ \ for all \ $j \in [0,3]$ \ and then we have
$$ a(1) = 4.x^{\alpha} +  \lambda.x^{\bf{1}} \quad {\rm and} \quad  b(1) = \vert\alpha\vert.x^{\alpha} + \lambda.x^{\bf{1}} $$
and then \ $(a - b)(1) = ( 4 - \vert\alpha\vert).x^{\alpha} $.\\
Using the other equalities gives after some easy computations
\begin{align*}
&   (a - 2b)(x^{\alpha}) =  ( 4 - \vert\alpha\vert).x^{\alpha}.\sigma(x^{\alpha}) \\
&   (a - 3b)\big(x^{\alpha}.\sigma(x^{\alpha})\big) =  ( 4 - \vert\alpha\vert).x^{\alpha}.\sigma(x^{\alpha}).\sigma^{2}(x^{\alpha}) \\
&   (a - 4b)\Big(x^{\alpha}.\sigma(x^{\alpha}).\sigma^{2}(x^{\alpha}))\Big) =  ( 4 - \vert\alpha\vert).x^{\alpha}.\sigma(x^{\alpha}).\sigma^{2}(x^{\alpha}).\sigma^{3}(x^{\alpha}).
\end{align*}
So we have
\begin{equation*}
(a - 4b)(a - 3b)(a - 2b)(a - b)[1] = ( \vert\alpha\vert - 4)^{4}.x^{\alpha}.\sigma(x^{\alpha}).\sigma^{2}(x^{\alpha}).\sigma^{3}(x^{\alpha}). \tag{@}
\end{equation*}

Now we compute the action of \ $b$ \ on the monomial \ $ m_{p} : = (x_{0}.x_{1}.x_{2}.x_{3})^{p}$. We have
$$B(m_{p}) = (p+1).b(m_{p}).v_{0} = M_{0}.m_{p}.\Sigma(x^{\alpha}) + \lambda.m_{p+1}.v_{0} $$
and so \ $m_{p}.\sigma^{j}(x^{\alpha}) $ \ is independent on \ $j \in [0,3]$ \ and equal to 
 $$(\vert\alpha\vert)^{-1}.\big((p+1).b(m_{p}) - \lambda.m_{p+1}\big).$$
Now we obtain
$$ \big(a - \frac{4(p+1)}{ \vert\alpha\vert }.b\big)(m_{p}) = \lambda.(\vert\alpha\vert - 4).m_{p+1} .$$
So we conclude that the element\footnote{The product is taken left to right for decreasing \ $p$.}
$$(a - 4b).\Big[\prod_{p=0}^{\vert\alpha\vert -2} \ \big(a - \frac{4(p+1)}{ \vert\alpha\vert }.b\big) - \lambda^{\vert\alpha\vert }.(\vert\alpha\vert - 4)^{\vert\alpha\vert - 4}.(a - 3b)(a - 2b)(a - b)\Big] $$
annihilated the class of \ $dx$ \ in \ $E_{f}$.\\

Note that this element depends on \ $\alpha$ \ only by the number \ $\vert\alpha\vert$. So the differential equation for the periods of   \ $[dx]$ \ for such an \ $f$ \ depends only on the total degree of the monomial \ $x^{\alpha}$ \ and not of the precise choice of such a monomial \ $x^{\alpha}$. \\
Remark also that for \ $\vert \alpha\vert = 5$ \ the element of \ $\Ab$ \ under the brackets is invariant by the anti-automorphism  \ $\theta_{4}$ \ of \ $\Ab$ \ defined by the relations
$$ \theta_{4}(x.y) = \theta_{4}(y).\theta_{4}(x),\quad \theta_{4}(1) = 1,\quad  \theta_{4}(a) = a - 4b,\quad  \theta_{4}(b) = -b .$$

So it is remarkable that we find the same differential equation with this symetry property for the period of \ $dx\wedge dy\wedge dz\wedge dt$ \ for the following polynomials :
\begin{align*}
&  x^{5} + y^{5} + z^{5} + t^{5} + \lambda.x.y.z.t \\
&  x^{4}.y + y^{4}.z + z^{4}.t + t^{4}.x + \lambda.x.y.z.t \\
&  x^{4}.z + y^{4}.t + z^{4}.x + t^{4}.y + \lambda.x.y.z.t \\
&  x^{3}.y.z + y^{3}.z.t +  z^{3}.t.x + t^{3}.x.y +  \lambda.x.y.z.t \\
&   x^{3}.y.t + y^{3}.z.x +  z^{3}.t.y + t^{3}.x.z +  \lambda.x.y.z.t \\
&   x^{2}.y^{2}.z + y^{2}.z^{2}.t + z^{2}.t^{2}.x + t^{2}.x^{2}.y  + \lambda.x.y.z.t \\
&   x^{2}.y^{2}.t + y^{2}.z^{2}.x + z^{2}.t^{2}.y + t^{2}.x^{2}.z  +  \lambda.x.y.z.t \\
&   x^{2}.z^{2}.t + y^{2}.t^{2}.x + z^{2}.x^{2}.y + t^{2}.y^{2}.z  + \lambda.x.y.z.t \\
& etc \dots
\end{align*}

\newpage

\section{The factorization theorems.}

\subsection{The decomposition theorem.}

Here we shall work with the completion \ $\A$ \ of the algebra \ $\Ab$ \ relative to its \ $b-$adic filtration given by the two-sided ideals \ $b^{\nu}.\Ab$ \ for \ $\nu \in \mathbb{N}$. \\
Recall that a (a,b)-module is a left \ $\A-$module which is free and of finite rank on the commutative sub-algebra \ $\C[[b]]$ \ of \ $\A$. For more on (a,b)-modules (generalized Brieskorn modules) see [B.06]. \\
We shall say that a (a,b)-module is monogenic when it may be generated by one element as a left \ $\A-$module.\\
Recall also that a (a,b)-module \ $E$ \ is {\em local}, respectively {\em simple pole}, if there exists an integer \ $N \geq 1$, respectively \ $N = 1$, such that \ $a^{N}.E \subset b.E$. \\

The purpose of this section is to prove a decomposition theorem for any (a,b)-module \ $E$. In the case where \ $E : = \A\big/\A.P$ \ where \ $P\in \Ab$ \ is monic in \ $a$, as in the differential equations constructed in the previous section, it will give a factorisation of  \ $P$ \ in \ $\A$
$$ P = P_{1}. \dots P_{k} $$
where each \ $P_{i}$ \ is monic in \ $a$ \ and with a class modulo \ $b.\A$ \ of the form \ $(a - v_{i})^{d_{i}}$ \ where \ $v_{i}$ \ is a complex number and \ $d_{i}$ \ a positive integer. So, in the decomposition theorem which corresponds to the ``decomposition'' of the differential equation in its different singular points, we find in this case local factors of the same type (i.e of the form \ $\A\big/\A.P_{i}$). \\

For a proof of the following proposition see [B.09].

 \begin{prop}\label{ann.}
  Let \ $I$ \ be a left ideal in \ $\A$ \ such that \ $E : = \A\big/I$ \ is a  rank \ $d$ \ (a,b)-module and  such that \ $x^d$ \ is the minimal polynomial of the action of \ $a$ \ on \ $E\big/b.E$. Then there exists a monic degree \ $d$ \ polynomial in \ $a$ \   with coefficients in \ $\C[[b]]$ \ such that \ $I = \A.P$. Moreover, such a \ $P$ \ is unique, and \ $E$ \ is regular if and only if the homogeneous  initial part of \ $P$ \ in (a,b) is of degree \ $d$.\\
   When \ $E$ \ is regular, the Bernstein element\footnote{The Bernstein element \ $ Q$ \ of \ $E$ \ is the homogeneous element in \ $\A$ \ monic in \ $a$ \ of degree \ $d$, such that the Bernstein polynomial of \ $E$ \ is given by \ $(-b)^{d}.B(-b^{-1}.a) = Q$.} of \ $E$ \ is the initial part of \ $P$.\\
  Conversely, if \ $E$ \ is a monogenic local (a,b)-module, for any generator \ $e$ \ of \ $E$, its annihilator in \ $\A$ \ is of the type \ $I = \A.P$ \ where \ $P$ \ is a monic polynomial in \ $a$ \ of degree \ $d$ \ with coefficients in \ $\C[[b]]$.
  \end{prop}
  
\bigskip

\begin{thm}[Decomposition] \label{decomp.}
Let \ $E$ \ be an (a,b)-module. For each spectral subspace \ $F_v$ \ for the eigenvalue \ $v$ \ of the action of \ $a$ \ on \ $F : = E\big/b.E$, there exists a normal\footnote{A sub-module \ $G \subset E$ \ is normal when it satisfies \ $G \cap b.E = b.G$.} submodule \ $G_v$ \ of \ $E$ \ such that \ $G_v\big/b.G_v \simeq F_v$, and such that \ $(a - v)^{d(v)}.G_v \subset b.G_v$, where \ $d(v)$ \ is the multiplicity of the root \ $v$ \ in the minimal polynomial of the action of \ $a$ \ on \ $F$. Moreover, \ $G_v$ \ is of rank \ $\dim_{\C} F_v$ \ and we have the  following decompositon  as a direct sum of  (a,b)-modules
$$ E =  \oplus_v \  G_v  .$$
\end{thm}

\parag{Proof} First we shall construct \ $G_v$ \ for a given eigenvalue \ $v$ \ of \ $a$ \ acting on \ $F$. We may assume without loss of generality that \ $v = 0$, and we note \ $d(0) = d$ \ for short. Define
$$ G : = \{ x \in E \ / \ a^{n.d}.x \in b^n.E \quad \forall n \in \mathbb{N}\} .$$
This vector subspace is clearly stable by \ $a$. It is also stable by \ $b$ \ using the following identities which are easy to prove by induction on \ $\nu \in \mathbb{N}$
\begin{align*}
& a^{\nu}.b = b.(a + b)^{\nu} \\
& (a + b)^{\nu} = a^{\nu} + \nu.a^{\nu-1}.b \\
& a^{\nu}.b = b.a^{\nu} + \nu.b.a^{\nu-1}.b 
\end{align*}
If \ $a^{(n+1).d}.x \in b^{n+1}.E$ \ and \ $a^{n.d}.b.x \in b^{n+1}.E$ \ then
$$ a^{(n+1).d}.b.x = b.a^{(n+1).d}.x + (n+1).d.b.a^{d-1}.(a^{n.d}.b.x) \in b^{n+2}.E .$$
So we obtain that \ $b.G \subset H$ \ where \ $H \subset G$ \ is defined by
$$ H : = \{ y \in E \ / \ a^{n.d}.y \in b^{n+1}.E \}.$$

We shall show now that for any \ $x \in E$ \ such that \ $x$ \ lies in \ $F_0 + b.E$, there exists an  element \ $\xi \in G$ \ such that \ $\xi - x \in b.E$.\\
Fix a vector space decomposition \ $E : = \tilde{F} \oplus b.E$ \ and denote \ $\tilde{F}_v$ \ the pull back of \ $F_v$ \ by the bijection \ $\tilde{F} \to F$. Define then
$$L : = \oplus_{v \not= 0}\tilde{F}_v .$$
 Assume we have already found \ $y_1, \dots, y_{n-1}$ \ in \ $L$ \ such that
\begin{equation*}
a^{n.d}.(x - \sum_{j=1}^{n-1} \ b^j.y_j ) = b^n.x_n \tag{$@_n$}
\end{equation*}
for some \ $x_n \in E$. For \ $n = 1$ \ this is clear from our assumption that \ $x \in F_0 + b.E$ \ and \ $a^d.F_0 \subset b.E$. Assume that \ $(@_n)$ \  is true for \ $n \geq 1$. We shall prove \ $(@_{n+1})$. Write
$$ x_n = t_n + z_n + b.w_n $$
where \ $t_n$ \ is in \ $\tilde{F}_0$, $z_n $ \ in \ $L$ \ and \ $w_n \in E$. As \ $a $ \ is bijective on \ $L = L + b.E\big/b.E $ \ it is also bijective on \ $b^{n}.L = b^n.L+ b^{n+1}.E\big/b^{n+1}.E$, because the identity 
$$a.b^n = b^n.a + nb^{n+1}$$ 
 gives the commutativity of the diagram
$$ \xymatrix{ E\big/b.E \ar[r]^a \ar[d]^{b^n} & E\big/b.E \ar[d]^{b^n} \\ b^n.E\big/b^{n+1}.E \ar[r]^a & b^n.E\big/b^{n+1}.E }$$
So we may find   \ $y_n \in \tilde{F}_v$ \ such that \ $a^{n.d}.b^n.y_n = b^n.(z_n + b.\zeta_n)$ \ where \ $\zeta_n $ \ is in \ $E$. Then we obtain
$$ a^{(n+1).d}.b^n.y_n = a^d.b^n.(z_n+ b.\zeta_n) $$
and so
\begin{align*}
&  a^{(n+1).d}.(x - \sum_{j=1}^n b^j.y_j) = a^d.b^n.x_n - a^{(n+1).d}.b^n.y_n  \\
& \qquad  = a^d.b^n.t_n + a^d.b^n.z_n  + a^d.b^{n+1}.w_n - a^d.b^n.(z_n + b.\zeta_n).
\end{align*} 
But \ $a^d.b^n.t_n = b^n.a^d.t_n + b^{n+1}.E$ \ and it belongs to  \ $b^{n+1}.E$ \ as \ $t_n \in \tilde{F}_0$ \ implies that \ $a^d.t_n$ \ is in \ $b.E$. So we conclude that  \ $a^{(n+1).d}.(x - \sum_{j=1}^n b^j.y_j)$ \ is in \ $b^{n+1}.E$. This proves the existence  of \ $ y : = \sum_{j=1}^{+\infty} b^j.y_j \in b.\C[[b]].L$ \ such that for any \ $n \in \mathbb{N}$ \ we have  $$a^{n.d}.(x - y) =  a^{n.d}.(x - \sum_{j=0}^{n-1} b^j.y_j) + a^{n.d}.b^n.\sum_{j=1}^{\infty} b^j.y_{n+j}  \in b^n.E .$$
So we have found an element \ $z : = x - y \in G$ \ such that \ $x - z \in b.E$. Note that this implies that the map \ $G\big/b.G \to E\big/b.E$ \ induced by the inclusion \ $G \subset E$ \  has exactely \ $F_0$ \ as image, because \ $a^d.x \in b.E$ \ implies that \ $x \in F_0 + b.E$.\\

Let us show now that \ $H = b.G$. As the inclusion \ $b.G \subset H$ \ is already proved, take  
\ $x \in H$. We have \ $x = b.z$ \ for some \ $z \in E$ \ (using \ $n=0$ \ in the definition of \ $H$), and now
$$ a^{n.d}.b.z = b.a^{n.d}.z + n.d.b.a^{n.d-1}.b.z = b.a^{n.d}.z + n.d.b.a^{d-1}.a^{(n-1).d}.b.z $$
and as \ $a^{(n-1).d}.b.z \in b^n.E$ \ and \ $a^{n.d}.b.z \in b^{n+1}.E$ \  (\ $x \in H$) \ we get
$$ b.a^{n.d}.z \in b^{n+1}.E \quad  {\rm and \ so} \quad a^{n.d}.z \in b^n.E .$$
 So \ $z $ \ is  in \ $G$. This implies that \ $a^d.G \subset H = b.G $.\\
 We shall prove now that \ $G$ \ is normal in \ $E$. We shall first show that if \ $x $ \ lies in \ $b^k.E\cap G$ \ for some \ $k \geq 1$ \ there exists\ $y \in G$ \ such that \ $x - b^k.y \in b^{k+1}.E$.\\
 Assume that \ $x \in G$ \ satisfies  \ $x = b^k.z$ \ for some \ $z \in E$ \ and \ $k \geq 1$. Write \ $z = u + v  + b.\xi $ \ with \ $u \in \tilde{F}_0, v\in L$ \ and \ $\xi \in E$. Then \ $x = b^k.z = b^k.u + b^k.v + b^{k+1}.\xi$ \ is the spectral decomposition relative to the action of \ $a$ \  on \ $x$ \ in \ $b^k.E\big/b^{k+1}.E$. This implies that \ $v = 0$ \ and so \ $x = b^k.u + b^{k+1}.\xi$ \ and the previous step of our proof gives an element \ $y \in G$ \ such that \ $y - u \in b.E$. So we have \ $x - b^k.y \in b^{k+1}.E$ \ proving our assertion.\\
 So, for such an \ $x \in b^k.E \cap G$ \ we can construct inductively \ $y_k, \dots, y_n, \dots$ \ in \ $G$ \ such that \ $x = b^k.\sum_{j=0}^{+\infty} b^j.y_{k+j} $. So \ $x $ \ is in \ $b^k.G$, and \ $G$ \ is normal.\\
 As a  consequence, the map \ $G\big/b.G \rightarrow E\big/b.E$ \ is  injective and is an isomorphism of \ $G\big/b.G$ \ on \ $F_0$.  So \ $G$ \ has rank \ $\dim_{\C}(F_0)$, and satisfies \ $a^d.G \subset b.G$.\\
 To prove that the sum of the normal sub-modules \ $G_{v}$ \ is direct is an easy exercice left to the reader. Then to prove the direct sum decomposition of \ $E$ \ it enough to prove that \ $E\big/b.E$ \ is the direct sum of the images of the \ $G_v\big/b.G_v$. But this is the spectral decomposition of \ $E\big/b.E$ \ for the action induced by \ $a$ \ from what we proved above.  $\hfill \blacksquare $\\
 
 \parag{Remark} As each \ $G_{v}$ \ is a direct factor (as (a,b)-module) of \ $E$, each of these (a,b)-modules is monogenic when \ $E$ \ is monogenic. $\hfill \square$\\
 
 The following corollary is now immediate.
 
 \begin{cor} \label{mono.}
 Let \ $P \in \A$ \ be a monic polynomial in \ $a$ \ with coefficients in \ $\C[[b]]$, and define \ $E : = \A\big/\A.P$. Let \ $v_{1}, \dots,v_{l}$ \ be the distinct  eigenvalues of \ $a$ \ acting on \ $E\big/b.E$, and let \ $d_{1}, \dots, d_{l}$ \ be the dimension of the corresponding spectral subspaces. Then there exists  monic polynomials \ $P_{1}, \dots, P_{l}$ \ respectiveley in \ $a - v_{1}, \dots, a - v_{l}$, of respective degrees \ $d_{1}, \dots, d_{l}$, with coefficients in \ $\C[[b]]$, such that we have in \ $\Ab$ \ the equality
 $$ P = P_{1}\dots P_{l}  .$$
 \end{cor}
 
 Of course, if we decompose the generator \ $e$ \ of \ $E$ \ according to the direct decomposition \ $E = \oplus_{v_{i}, i \in [1,l]} G_{v_{i}}$ :
 $$ e = \sum_{i=1}^{l} \ e_{v_{i}} $$
 then \ $P_{l}$ \ generates the annihilator of \ $e_{v_{l}}$ \ in \ $E$ \ (or in \ $G_{v_{l}}$). So we have an isomorphism  \ $\A\big/\A.P_{l} \simeq G_{v_{l}}$ \ by sending \ $1$ \ to \ $e_{v_{l}}$.\\
  Again the polynomial \ $P_{l-1}$ \ generates the annihilator of \ $P_{l}.e_{v_{l-1}}$. As \ $P_{l}.e_{v_{l-1}}$ \ is a generator of \ $G_{v_{l-1}}$ \ because in a local (a,b)-module the topology defined by  the \ $a-$filtration is complete, so the action of an element in \ $\A$ \ with a non zero constant term is invertible, we have an isomorphism  \ $\A\big/\A.P_{l-1} \simeq G_{v_{l-1}}$ \  obtained by sending \ $1$ \ to \ $P_{l}.e_{v_{l-1}}$ \ etc... 
 $\hfill \blacksquare$\\
 
   Remark that changing the order of the \ $v_{i}$ \ changes the polynomials \ $P_{i}$ \ in the decomposition above but only change the order of the set of the isomorphism classes of the \ $\A\big/\A.P_{i}$ .\\
 
 The following consequence of the theorem \ref{base} and the results above is  obvious.

 \begin{cor}\label{manquant}
 In the situation of the theorem \ref{base} the element \ $P_{d}$ \ is a left multiple of the Bernstein element of the left \ $\A-$module generated by \ $[\mu.dx]$ \ in \ $E_{f}\otimes_{\C[b]}\otimes \C[[b]]\big/b-torsion$.
 \end{cor}
 
 This is, of course, a non trivial information on the monodromies at the origin for the periods integrals considered in the beginning of this paper.

\subsection{Irregularity and the second factorization theorem.}

In this paragraph we examine what happens in the previous factorization when the initial part of \ $P$ \ is not monic in \ $a$.

\begin{defn}\label{Irr.0}
A (a,b)-module \ $E$ \ will be called  {\bf totally irregular (at \ $0$)} if any \ $\A-$linear map \ $f : E \to F$ \ in a simple pole (a,b)-module \ $F$ \ is the zero map.
\end{defn}

This terminology is compatible with the usual notion of irregularity of a differential equation at the point \ $0$. It concerns in fact only \ $G_{0}$ \ the local part at \ $0$ \ in the previous decomposition theorem \ref{decomp.} for an (a,b)-module.\\

A quotient of a totally irregular (a,b)-module is again totally irregular.

\begin{lemma}\label{Irr.1}
Let \ $E$ \ be a local (a,b)-module. Then there exists a smallest normal sub-module \ $I \subset E$ \ such that the quotient \ $E\big/I$ \ is regular. This sub-module \ $I$ \  is  totally irregular.
\end{lemma}

\parag{Proof}  First remark that if \ $I$ \ and \ $J$ \ are normal sub-modules such that \ $E\big/I$ \ and \ $E\big/J$ \ are regular, then \ $E\big/I\cap J$ \ is again regular because \ $E\big/I \cap J$ \ is a sub-module of the direct sum of the regular modules \ $E\big/I$ \ and \ $E\big/J$. Then let  \ $I$ \ be the intersection of all normal  sub-modules \ $J$ \  of \ $E$ \ such that  \ $E\big/J$ \ is regular. As any descending chain of normal sub-modules is finite, it is clear that \ $I$ \ is the smallest normal sub-module of \ $E$ \ such that  \ $E\big/I$ \ is regular. Let \ $ f : I \to F$ \ be a \ $\A-$linear map, where \ $F$ \ has a simple pole, and let \ $K$ \ be the kernel of \ $f$. Then in the exact sequence of (a,b)-modules
$$0 \to I\big/K \to E\big/K \to E\big/I \to 0 $$
the (a,b)-modules \ $I\big/K$\footnote{As \ $f$ \ induces an injection of \ $I\big/K$ \ in \ $F$ \ which has a simple pole, it is regular.} \ and \ $E\big/I$ \ are regular. So is \ $E\big/K$. Then \ $E\big/K$ \ is a sub-module of a simple pole (a,b)-module, and then \ $I \subset K$. This implies \ $f = 0$, and so \ $I$ \ is totally irregular. $\hfill \blacksquare$

\begin{lemma}\label{Irr.2}
Let \ $P \in \A$ \ be a monic polynomial in \ $a$ \ of degree \ $d$. Assume that the initial form of \ $P$ \ in (a,b) is equal to \ $b^{q}$, for some integer \ $q < d $.
 Then the (a,b)-module \ $E : = \A \big/ \A.P$ \ is totally irregular.
 \end{lemma}

\parag{Proof} Let \ $f : E \to F$ \ be a \ $\A-$linear map such that \ $F$ \ is a simple pole (a,b)-module. Let \ $x : = f(1)  \in F$. We have \ $P.x = 0$ \ in \ $F$, and so \ $b^{q}.x $ \ is in \ $b^{q+1}.F$, because of our assumption on the initial form of \ $P$ \ and on the simple pole for \ $F$. So \ $ x \in b.F$. Then the image of \ $f$ \ lies in \ $b.F$ \ which is again a simple pole (a,b)-module. Iteration of this gives that the image of \ $f$ \ is in  \ $ b^{n}.F$ \ for any integer \ $n$. But \ $\cap_{n \in \mathbb{N}} \ b^{n}.F = 0$, and so \ $f = 0$.$\hfill \blacksquare $\\

\begin{thm}\label{Irr.3}
Let \ $P$ \ be a element in \ $\A$ \ with the following properties :
\begin{enumerate}[(i)]
\item \ $P$ \ is a monic polynomial in \ $a$ \ of degree \ $d+h$ \ with \ $h \geq 1$.
\item The initial form in (a,b) of \ $P$ \ has degree \ $d \geq q \geq 0$, and is of the form \ $\rho.b^{q}.P_{d-q}$ \ with \ $ \rho \in \C^{*}$ \ where \ $P_{d-q}$ \ is monic in \ $a$.
\item Modulo \ $b.\A$, we have \ $P = a^{d+h}$.
%\item Assume that the polynomial \ $B$ \  associated to \ $P_{d-q}$ \ by the relation \ $b^{d-q}.B(-b^{-1}.a) = P_{d-q}$ \ in \ $\Ab$, has negative rational roots.
\end{enumerate}
Then there exists elements \ $Z$ \ and \ $Q$ \ in \ $\A$ \ which are polynomials in \ $a$ \ of respective degrees \ $q+h$ \ and at most \ $ d-q-1$, with respective valuations in (a,b) at least \ $q+1$ \ and \ $d-q+1$, such that we have in \ $\A$ \ the equality
$$( \rho.b^{q} + Z).(P_{d-q} + Q) = P .$$
Defining the (a,b) modules 
 $$E : = \A\big/\A.P, \quad F : = \A\big/\A.( P_{d-q} + Q) \quad \rm{and}  \quad  F  : = \A\big/\A.( \rho.b^{q} + Z),$$
 we have the exact sequence of (a,b)-modules of respective ranks \ $q+h, d+h, d-q$ :
$$ 0 \to I \to E \to F \to 0 $$
and \ $I$ \ is the totally irregular part of \ $E$.
\end{thm} 

\parag{Remark} Assuming that the polynomial \ $B$ \  associated to \ $P_{d-q}$ \ by the relation \ $(-b)^{d-q}.B(-b^{-1}.a) = P_{d-q}$ \ in \ $\Ab$, has negative rational roots, the (a,b)-module \ $F$ \ is a geometric (a,b)-module generated by one element as left  \ $\A-$module ( see [B.09]) and \ $B$ \ is the Bernstein polynomial of \ $F$.$\hfill \square$

\parag{Proof} We shall construct by induction on \ $ n \geq d+1$ \ homogeneous elements in (a,b) \ $ \xi_{n-d+q}$ \ and \ $\eta_{n-q}$ \ of respective degrees \ $n-d+q$ \ and \ $n-q$, which are polynomial in \ $a$ \ of respective degrees less or equal to \ $q+h$ \ and \ $d-q-1$, such that for each \ $N \geq d+1$ \ we have
\begin{equation*} 
(b^{q} + \sum_{n =d+1}^{N} \ \xi_{n-d+q} ).(P_{d-q} + \sum_{n =d+1}^{N} \ \eta_{n-q}) = P + Y_{N+1} \tag{$@_{N}$}
\end{equation*} 
where \ $Y_{N+1}$ \ has valuation in (a,b)  \ $\geq N+1$ \ and is a polynomial in \ $a$ \ of degree \ $\leq d+h$.\\
For \ $N = d+1$ \ denote \ $(P )_{d+1}$ \ the homogeneous part in (a,b) of degree \ $d+1$ \ of \ $P$, and make the right division by \ $P_{d-q}$ \ in \ $\Ab$ :
$$ (P)_{d+1} = \xi_{q+1}.P_{d-q} + R_{d+1} $$
where \ $R_{d+1}$ \ has degree \ $\leq d-q-1$ \ in \ $a$ \ and is homogeneous in (a,b) of degree \ $d+1$. So we may write \ $ R_{d+1} =  \rho.b^{q}.\eta_{d-q+1}$. 
Let \ $Y_{d+2}^{0} : = P - \rho.b^{q}.P_{d-q} - (P )_{d+1}$. Then we have
\begin{align*}
& (\rho.b^{q} + \xi_{q+1}).(P_{d-q} + \eta_{d-q+1}) = \rho.b^{q}.P_{d-q} + \rho.b^{q}.\eta_{d-q+1} + \xi_{q+1}.P_{d-q} +  \xi_{q+1}.\eta_{d-q+1} \\
& \qquad = \rho.b^{q}.P_{d-q} + (P)_{d+1} + \xi_{q+1}.\eta_{d-q+1} = P + Y_{d+2}^{0} +  \xi_{q+1}.\eta_{d-q+1}  = P + Y_{d+2}
\end{align*}
where \ $Y_{d+2}$ \ has a valuation in (a,b) \ $ \geq d+2$ \ and degree \ $\leq d+h$ \ in \ $a$.\\
Assume now that we have proved \ $(@_{N})$ \ for some \ $N \geq d+1$. Let \ $Y_{N+1}^{0}$ \ the homogeneous part of \ $Y_{N+1}$ \ of degree \ $N+1$ \ in (a,b), and write \ $Y_{N+2}^{0}: = Y_{N+1}- Y_{N+1}^{0}$. The right division by \ $P_{d-q}$ \ in \ $\A$ \ gives 
$$ -Y_{N+1}^{0} = \xi_{N+1-d+q}.P_{d-q} + R_{N+1} $$
where \ $\xi_{N+1-d+q}$ \ is homogeneous in (a,b) of degree \ $N+1-d+q$, where \ $R_{N+1}$ \ is homogeneous in (a,b) of degree \ $N+1$, and they have respectively degrees less than \ $q+h$ \ and \ $d-q-1$ \ in \ $a$. Then we may write \ $R_{N+1} = \rho.b^{q}.\eta_{N-q+1}$. Now we have
\begin{align*}
& (\rho.b^{q} +  \sum_{n =d+1}^{N+1} \ \xi_{n-d+q} ).(P_{d-q} + \sum_{n =d+1}^{N+1} \ \eta_{n-q})  \\
& \qquad = P + Y_{N+1} +  \xi_{N+1-d+q}.P_{d-q} + \rho.b^{q}.\eta_{N-q+1} + Z_{N+2} \\
& \qquad = P + Y_{N+2}^{0} + Z_{N+2} = P + Y_{N+2} 
\end{align*}
where we defined 
 $$Z_{N+2} : = \xi_{n-d+q} ).\sum_{d+1}^{N+1} \eta_{n-q} + (\sum_{d+1}^{N} \xi_{n-d+q}).\eta_{N+1-q}$$ 
  and \ $Y_{N+2} : = Y_{N+2}^{0} + Z_{N+2}$. This complete the proof of our induction step.\\
Now the series \ $ \sum_{n = d+1}^{\infty} \ \xi_{n-d+q}$ \ and \ $\sum_{n = d+1}^{\infty} \ \eta_{n-q}$ \ converge in \ $\A$ \ to  \ $Z$ \ and \ $Q$ \ which are polynomials in \ $a$ \ of degree\footnote{Remark that \ $\xi_{q+h}$ \ has degree exactly \ $q+h$ \ and that \ $\xi_{m}$ \ has degree \ $\leq q+h-1$ \ for \ $m \not= q+h$.} respectively equal to \ $q+h$ \ and \ $\leq d-q-1$. So we obtain the desired  decomposition of \ $P$ \ 
by defining \ $Z : =  \sum_{n = d+1}^{\infty} \ \xi_{n-d+q}$ \ and \ $Q : = \sum_{n = d+1}^{\infty} \ \eta_{n-q}$.\\
Now the lemma \ref{Irr.2} implies that \ $I$ \ is totally irregular and the lemma \ref{Irr.1} implies that \ $E\big/I$ \ is regular. To complete the  proof is then easy. $\hfill \blacksquare$\\

\parag{Example} Let \ $P_{d+h}$ \ and \ $P_{d-q}$ \ be homogeneous elements in \ $\Ab$ \ with respective degree \ $d+h$ \ and \ $d$. Assume that \ $P_{d+h}$ \ and \ $P_{d-q}$ \ are monic in \ $a$, then \ $P : = P_{d+h} + \rho.b^{q}.P_{d-q}$ \ satisfies the hypothesis of the theorem.$\hfill \square$\\

In the case of the previous example, let \ $P_{d+h} = \xi_{h+q}.P_{d-q} + \rho.b^{h+q+1}.\tilde{\eta}_{d-q-1}$ \ be the right division of \ $P_{d+h}$ \ by \ $P_{d-q}$ \ in \ $\Ab$. Then we have \ $ Z = \rho.b^{q} + \xi_{h+q} + Y$ \ with \ $Y$ \ a degree \ $\leq h+q-1$ \ polynomial in \ $a$ \ with valuation in (a,b) at least equal to \ $2h+q$, and \ $Q = P_{d-q} + b^{h+1}.\tilde{\eta}_{d-q-1} + T$ \ where \ $\tilde{\eta}_{d-q-1}$ \ and  \ $T$ \ are  degree \ $\leq d-q-1$ \ polynomials in \ $a$ \ with valuation in (a,b) respectively equal to \ $d-q-1$ \ and at least equal to \ $d+2h-q$.

\parag{Consequence} If an element  \ $x$ \ of an  (a,b)-module \ $E$ \  is killed by an element in \ $\A$ \ which is monic in \ $a$ \ of degree \ $d+h, h \geq 1$, equal to \ $a^{d+h}$ \ modulo \ $b.\A$,  with initial form equal to \ $\rho.b^{q}.P_{d-q}$, where \ $\rho$ \ is in \ $ \C^{*}$, where \ $q \geq 1$ \ and where \ $P_{d-q}$ \ is homogeneous of degree \ $d-q$ \ in (a,b) and monic in \ $a$, then the left \ $\A-$module  generated by \ $x$ \ in \ $E$ \ is  a regular (a,b)-module with rank at most \ $d-q$ \ and its Bernstein element is a right divisor of \ $P_{d-q}$ \ in \ $\Ab$. $\hfill \square$

\parag{References}

\begin{itemize}

\item{[B.06]} Barlet, D. {\it Sur certaines singularit\'es non isol\'ees d'hypersurfaces I}, Bull. Soc. math. France 134 (2), ( 2006), p.173-200.

%\item{[B.II]} Barlet, D. {\it Sur certaines singularit\'es d'hypersurfaces II}, J. Alg. Geom. 17 (2008), p. 199-254.

\item{[B.09]} Barlet,D. {\it P\'eriodes \'evanescentes et (a,b)-modules monog\`enes}, Bollettino U.M.I. (9) II (2009), p. 651-697.

\item{[Br.70]} Brieskorn, E. {\it Die Monodromie der Isolierten Singularit{\"a}ten von Hyperfl{\"a}chen}, Manuscripta Math. 2 (1970), p. 103-161.

%\item{[D.70]} Deligne, P. {\it Equations diff\'erentielles \`a points singuliers r\'eguliers} Lecture Notes in Math. 163, Springer 1970.

%\item{[G.65]} Grothendieck, A. {\it On the de Rham cohomology of algebraic varieties} Publ. Math. IHES 29 (1966), p. 93-101.

% \item{[K.76]} Kashiwara, M. {\it b-function and holonomic systems}, Inv. Math. 38 (1976) p. 33-53.

\item{[M.74]} Malgrange, B. {\it Int\'egrale asymptotique et monodromie}, Ann. Sc. Ec. Norm. Sup. 7 (1974), p. 405-430.

\item{[M.75]} Malgrange, B. {\it Le polyn\^ome de Bernstein d'une singularit\'e isol\'ee}, in Lect. Notes in Math. 459, Springer (1975), p. 98-119.

\item{[Mi.68]  Milnor, J.  \textit{Singular Points of Complex Hypersurfaces .} Ann. of Math. Studies  61 (1968) Princeton  .}

\end{itemize}

\end{document}